# Frontier improvement in the DEA models

by

Vladimir E. Krivonozhko[a,b,*], Finn R. Førsund[c], Andrey V. Lychev[a]

[a] National University of Science and Technology MISiS, Leninskiy prs. 4, Moscow 119049, Russia

[b] Federal Research Center "Computer Science and Control" of the Russian Academy of Sciences, Vavilov st. 44-2, Moscow 119333, Russia

[c] Department of Economics, University of Oslo, 0317 Blindern Norway

**Abstract:** Applications of data envelopment analysis (DEA) show that many inefficient units are projected onto the weakly efficient parts of the frontier when efficiency scores are computed. However this fact disagrees with the main concept of the DEA approach, because the efficiency score of an inefficient unit has to be measured relative to an efficient unit. As a consequence inaccurate efficiency scores may be obtained. This happens because a non-countable (continuous) production possibility set is determined on a basis of a finite number of production units. It has been proposed in the literature to use artificial production units in the primal space of inputs and outputs as a starting point in order to improve the frontier of the DEA models. Farrell was the first who introduced artificial units in the primal space of inputs and outputs in order to secure convex isoquants. In previous papers we introduced the notion of terminal units. Moreover, some relationships were established between terminal units and other sets of units that were proposed for improving envelopment. In this paper we develop an algorithm for improving the frontier. The construction of algorithm is based on the notion of terminal units. Our theoretical results are verified by computational experiments using real-life data sets and also confirmed by graphical examples.

[*] Corresponding author. Tel. +74956384473, mail address: National University of Science and technology MISiS, Leninskiy prs. 4, Moscow 119049, Russia
E-mail addresses: krivonozhko@gmail.com (V.E. Krivonozhko), f.r.forsund@econ.uio.no (F.R. Førsund), lychev@misis.ru (A.V. Lychev).

# 1. Introduction

The estimation of efficiency scores using the non-parametric data envelopment analysis (DEA) model is based on comparing observations with points on the frontier. The term frontier is used because it is based on best practice observations.

However, the DEA frontier is an estimate of the real unobserved technology. Furthermore, the estimated DEA production possibility set is assumed to be a convex polyhedral set with the closest fit to the data. The efficient border of the set represents the frontier. However, a border formed by facets spanned by extreme efficient units (vertex points) is not based on empirical insight of the activity in question, but is a consequence of the method. The DEA frontier is based on using all observed production units. It is assumed to be common for all units.

It is rather common that detailed engineering information about production processes is represented by production functions of a rather smooth type facilitating estimation of interesting concepts for managers, like marginal productivities, rate of substitution between inputs, and rate of transformation between outputs. As Farrell (1957, p. 262) remarks: "It will not be surprising if the method of estimation is not the best for any particular use, for it was chosen simply as providing the best measure of technical efficiency." Improving the frontier is therefore based on trying to make the frontier more accommodating for the most common properties specified by analysts.

Applications of DEA show that many inefficient units are projected on the weakly efficient parts of the frontier when efficiency scores are computed (see Section 2 for definitions). However, this fact disagrees with the main concept of the DEA approach, because the efficiency score of an inefficient unit has to be measured relative to efficient units. As a consequence inaccurate efficiency scores may be obtained. This happens because a non-countable (continuous) production possibility set is determined on a basis of a finite number of production units.

The question is how to do improvements. Past efforts started with Farrell (1957) introducing artificial units of zero and infinity for each input dimension in his non-parametric constant returns to scale (CRS) model in order to ensure convexity of the piecewise linear isoquants. Thanassoulis and Allen (1998) introduced a set of extreme efficient points termed anchor units in the case of constant returns to scale and a single input, having previously explored the use of weight restrictions as a way of overcoming inadequacies with the DEA method (see e.g. Førsund (2013) for a review of weight restriction approaches).

Allen and Thanassoulis (2004) elaborated further the idea of focusing on anchor units as points of departure for formulating coordinates of artificial units. The main purpose was to improve the envelopment by reducing the number of inefficient units not "properly" enveloped, meaning projections to the frontier of these units being on a weakly efficient part of the frontier (see Section 2 for the formal definition).

Bougnol and Dulá (2009) introduced their definition of anchor units for the case of variable returns to scale and multiple inputs and outputs. They also proposed algorithm for finding anchor units. However, their algorithm may generate efficient units that are not anchor units in DEA models.

Thanassoulis et al. (2012) developed further the super-efficiency method for discovering anchor units in Banker, Charnes and Cooper (1984) model (hereafter called BCC). However, their method does not reveal all efficient units that may be the point of departure for improving the frontier in BCC models. Furthermore, their definition and their model produce different sets of anchor units. Moreover, their method of frontier improvement may turn initially efficient units into inefficient ones.

Edvardsen et al. (2008) developed an empirical method for determining some units that may generate incorrect results in the DEA models, they called them exterior units. At the same time, their methods cannot find all such units.

Krivonozhko et al. (2009) showed that incorporation of domination cones in the dual space of multipliers in DEA models for improving the frontier is equivalent to incorporation of artificial units and rays in the primal space of inputs and outputs which makes the process of frontier improving more visible and understandable.

Krivonozhko et al. (2015, 2017) defined the notion of terminal units. It was proved that the set of terminal units either contains all sets of units proposed in the scientific literature for improving the frontier or that some sets of units may contain redundant units that are not suitable for frontier improvement. On the other hand, every terminal unit can be used for improving the frontier. All these results are followed from the theorems in Krivonozhko et al. (2015, 2017). Moreover, some relationships between different sets of units (different sets of anchor units, exterior units and terminal units) that may cause biased results in DEA models were established.

A simple illustration may help to explain the background for introducing artificial units. In the two- dimensional Fig. 1 in input space, $L_2$-$L_1$ is the frontier isoquant for units *G* and *A*. Units *E*, *D* and *C* are the efficient units and the segments *E*-$L_2$ and *C*-$L_1$ are the weakly efficient faces. There is no support in the data for these two vertical and horizontal segments, respectively; they just represent the default routine of the estimation of the frontier of providing the most pessimistic estimate of the frontier technology given the assumption of convexity. This means that the efficiency score for unit *A*, *OQ/OA* where the projection point of unit *A* onto the frontier is *Q*, is the best estimate of efficiency the DEA model will give. However, the efficiency measure for unit G is *OR*/*OG*, and R is not on an efficient segment, but on a weakly efficient edge. This edge being the most pessimistic estimate of the frontier technology yields the efficiency measure as the most optimistic one. But since the segment $L_2E$ is an arbitrary choice without basis in data the efficiency measure may be way off from a true measure. If the true technology is expressed by the segment $L_3E$ then the efficiency score for unit G is reduced to *OS*/*OG*. So, unit S is inserted in the model as an artificial unit and unit E ceases to be terminal unit, but stays an efficient unit.

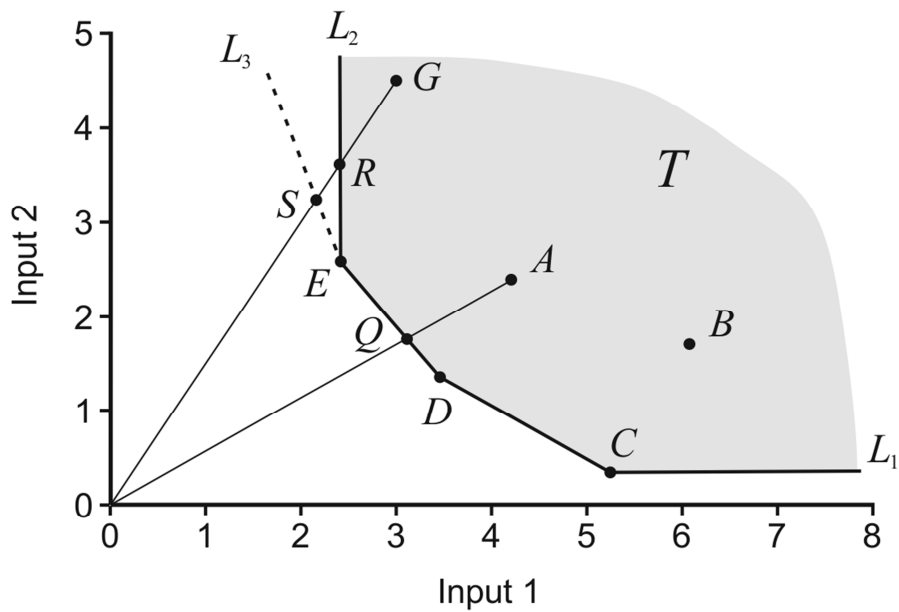

Figure 1. Two-dimensional input isoquant.

The purpose of this paper is to develop a general algorithm for improving the frontier in DEA models. We take the notion of terminal units as a point of departure for construction of the algorithm. Our theoretical results are verified by computational experiments using real-life data sets and also illustrated by graphical examples.

The paper is organized in the following way. In Section 2 the DEA models used and key definitions are stated, and in Section 3 how to use terminal units for constructing artificial units is outlines and illustrated. The algorithm is developed in Section 4. Numerical experiments are presented in Section 5 and illustrated based on sets of data from Russian banks and Norwegian municipal nursing and home care services. Section 6 concludes.

## 2. Background

In order for the paper to be self-contained we present the models we use and key definitions. It was shown in the DEA scientific literature (Krivonozhko et al., 2009) that the BCC model exhibiting variable returns to scale can approximate any DEA model from a large family of DEA models. For this reason, we consider the BCC model as a basic model in our exposition.

Consider a set of $n$ observations of actual production units $(X_j, Y_j)$, $j = 1,...,n$, where the vector of outputs $Y_j = (y_{1j},..., y_{rj}) \geq 0$, $j = 1,...,n$, is produced from the vector of inputs $X_j = (x_{1j},..., x_{mj}) \geq 0$. The production possibility set $T$ is the set $\{(X,Y) \mid \text{the outputs } Y \geq 0 \text{ can be produced from the inputs } X \geq 0\}$. The primal input-oriented BCC model can be written in the form

$$\min \theta$$

subject to

$$\sum_{j=1}^{n} X_j \lambda_j + S^- = \theta X_o,$$

$$\sum_{j=1}^{n} Y_j \lambda_j - S^+ = Y_o,$$

$$\sum_{j=1}^{n} \lambda_j = 1, \quad (1a)$$

$$\lambda_j \geq 0, \quad j = 1,...,n,$$

$$s_k^- \geq 0, \quad k = 1,...,m,$$

$$s_i^+ \geq 0, \quad i = 1,...,r,$$

where $X_j = (x_{1j},...,x_{mj})$ and $Y_j = (y_{1j},...,y_{rj})$ represent the observed inputs and outputs of production units $j = 1,...,n$, $S^- = (s_1^-,...,s_m^-)$ and $S^+ = (s_1^+,...,s_r^+)$ are vectors of slack variables. In this primal model the efficiency score $\theta$ of the production unit $(X_o, Y_o)$ is found; $(X_o, Y_o)$ is any unit from the set of production units $(X_j, Y_j)$, $j = 1,...,n$.

Notice that we do not use an infinitesimal constant $\varepsilon$ (a non-Archimedean quantity) explicitly in the DEA models, since we suppose that each model is solved in two stages in order to separate efficient and weakly efficient units.

The BCC primal output-oriented model can be written in the following form

$$\max \eta$$

subject to

$$\sum_{j=1}^{n} X_j \lambda_j + S^- = X_0,$$

$$\sum_{j=1}^{n} Y_j \lambda_j - S^+ = \eta Y_0,$$

$$\sum_{j=1}^{n} \lambda_j = 1, \quad (1b)$$

$$\lambda_j \geq 0, \quad j = 1,...,n,$$

$$s_k^- \geq 0, \quad k = 1,...,m,$$

$$s_i^+ \geq 0, \quad i = 1,...,r.$$

The production possibility set (PPS) for the BCC model can be written in the form

$$T = \left\{ (X,Y) \,\middle|\, \sum_{j=1}^{n} X_j \lambda_j \leq X, \sum_{j=1}^{n} Y_j \lambda_j \geq Y, \sum_{j=1}^{n} \lambda_j = 1, \lambda_j \geq 0, j = 1,\ldots,n \right\}. \quad (2)$$

The production possibility set $T$ in (2) is a convex polyhedral set. According to the classical theorems of Goldman (1956) and Motzkin (1936) any convex polyhedral set can be represented as a vector sum of convex combination of vertices and the non-negative linear combination of vectors (rays).

Any point of $T$ can be considered as a production unit in the process of decision-making and can be evaluated by model (1a) and model (1b). For this reason, in what follows we will use terms unit and point interchangeably.

Let $\theta^*$ be the optimal value of objective function in model (1a).

**Definition 1.** (Cooper et al. 2007). *Unit $(X_o, Y_o) \in T$ is called efficient with respect to the input-oriented BCC model if any optimal solution of (1a) satisfies: a) $\theta^* = 1$, b) all slacks $s_k^-$, $s_i^+$, $k = 1,\ldots,m$, $i = 1,\ldots,r$ are zero.*

All units satisfying condition (a) are referred as input weakly efficient with respect to the BCC input-oriented model. We denote the set of these weakly efficient points by $WEff_I T$. In the DEA literature (Seiford & Thrall, 1990; Banker & Thrall, 1992) this set is also called the input boundary.

Let $\eta^*$ be the optimal value of objective function in model (1b).

**Definition 2.** (Cooper et al. 2007). *Unit $(X_o, Y_o) \in T$ is called efficient with respect to the output-oriented BCC model if any optimal solution of (1b) satisfies: a) $\eta^* = 1$, b) all slacks $s_k^-$, $s_i^+$, $k = 1,\ldots,m$, $i = 1,\ldots,r$ are zero.*

All units satisfying condition (a) are referred as output weakly efficient with respect to the BCC output-oriented model. We denote the set of these weakly efficient points by $WEff_O T$. In the DEA literature (Seiford & Thrall, 1990; Banker & Thrall, 1992), this set is also called the output boundary.

**Definition 3.** (Banker & Thrall, 1992). *Activity $(X', Y') \in T$ is weakly Pareto efficient if and only if there is no $(X, Y) \in T$ such that $X < X'$ and $Y > Y'$.*

We denote the set of weakly Pareto efficient activities by $WEff_P T$, and the set of efficient points of $T$ with respect to the BCC model (1) by $Eff\ T$. Krivonozhko et al. (2005) have proved that the following relations hold:

$$Eff\ T \subseteq WEff_I T \cap WEff_O T,\ WEff_I T \cup WEff_O T \subseteq WEff_P T = Bound\ T,$$

where the boundary of T is designated as *Bound T*.

Before going further, let us recall some notions from convex analysis. Faces are formed by an intersection of the supporting hyperplane and the polyhedral set. In the DEA models, the dimension of face may vary from 0 up to $(m+r-1)$, the maximal dimension. Faces of maximal dimension are called facets. Faces of 0-dimension are known as vertices, 1-dimension as edges.

**Definition 4.** (Krivonozhko et al., 2015). *We call an extreme efficient unit terminal unit if an infinite edge is going out from this unit.*

According to Krivonozhko et al. (2015) only vectors of the following forms: $d_k = (e_k, 0) \in E^{m+r}$, $k=1,\ldots,m$, $g_i = -(0, e_i) \in E^{m+r}$, $i=1,\ldots,r$ can be the direction vectors of infinite edges of set $T$, where $e_k = (0,\ldots 1,\ldots 0) \in E^m$ (the unity is in *k*-th position) and $e_i = (0,\ldots,1,\ldots 0) \in E^r$ (the unity is in *i*-th position). A set of such direction vectors for given terminal unit we call *terminal directions* associated with this unit.

We denote the set of terminal units with respect to the production possibility set (2) by $T_{term}$. The models for determination of all terminal units of set $T$ are given in Krivonozhko et al. (2015). The relationship between $T_{term}$ and other sets of units proposed for improving the frontier is established in Krivonozhko et al. (2017), where the following assertion was formulated:

$$T_{ext} \subseteq T_{anc}^3 \subseteq T_{term} \subseteq T_{anc}^1,$$
$$T_{anc}^2 \subseteq T_{term} \subseteq T_{anc}^1.$$

In these relations $T_{anc}^1$ denotes the set of anchor units with respect to the definition of Bougnol and Dulá (2009), $T_{ext}$ denotes the set of exterior units (Edvardsen et al. 2008), $T_{anc}^2$ and $T_{anc}^3$ are the sets of anchor units according to the definition in Thanassoulis et al. (2012, p.178) and generated by their model, respectively.

## 3. Preliminaries

Under the development of the algorithm for improving the frontier we stick to the following principles:
a) all originally efficient units have to stay efficient after the frontier transformation;
b) every inefficient unit will be projected on to the efficient part of the frontier.

First of all, all terminal units are determined. Models for discovering of such units are described in Krivonozhko et al. (2015). Then two-dimensional sections are constructed for every terminal unit. For our purposes we need three types of sections.

Let us define a section of the frontier with a two-dimensional plane (Krivonozhko et al., 2004)

$$Sec(X_o, Y_o, d_1, d_2) = \{(X,Y) \mid (X,Y) \in Pl(X_o, Y_o, d_1, d_2) \cap WEff_P T\},$$

where $Pl(X_o, Y_o, d_1, d_2)$ is a two-dimensional plane going through point $(X_o, Y_o)$ and it is spanned by vectors $d_1, d_2 \in E^{m+r}$.

In our exposition we will use the following three types of sections.

1. Input isoquant, section $S_1$. In this case we take the following directions $d_1 = (e_i, 0) \in E^{m+r}$, $d_2 = (e_j, 0) \in E^{m+r}$, where $e_i$ and $e_j$ are $m$-identity vectors with a one in position $i$ and $j$, respectively. Hence this isoquant is determined by input variables $x_i$ and $x_j$.

2. Output isoquant, section $S_2$. In this case vectors for cutting the frontier are determined as follows: $d_1 = (0, e_i) \in E^{m+r}$, $d_2 = (0, e_j) \in E^{m+r}$, $e_i$ and $e_j$ are $r$-identity vectors with a one in position $i$ and $j$, respectively. This isoquant is associated with output variables $y_i$ and $y_j$.

3. Section $S_3$ reflects the dependence between variables $y_p$ and $x_s$. For construction of such dependence we took directions: $d_1 = (0, e_p) \in E^{m+r}$, where $e_p$ is $r$-identity vector with a one in position $p$, $d_2 = (e_s, 0) \in E^{m+r}$, $e_s$ is $m$-identity vector with a one in position $s$.

Figure 2 represents an input isoquant for some terminal unit $Z_k$, where $x_i$ and $x_j$ are input variables. If artificial unit $A$ is inserted somewhere in the region limited by rays $Z_k B$, $Z_k C$ and axis $Ox_i$, then unit $Z_k$ becomes just an efficient unit. Such operations can be accomplished for every terminal unit and for every type of sections going through this unit and that were described above. Observe that components of artificial unit $A$ coincide with corresponding components of unit $Z_k$, except coordinates that correspond to variables $x_i$ and $x_j$. In other words, unit $A$ belongs to the section that is going through point $Z_k$ and is determined by variables $x_i$ and $x_j$.

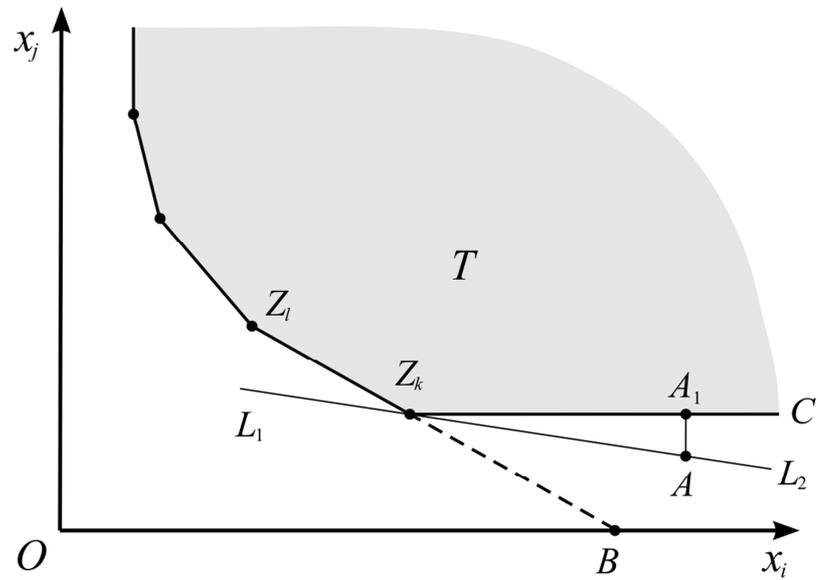

Figure 2. Terminal unit $Z_k$ turns into just an efficient unit

Figure 3 depicts an output isoquant. The two-dimensional section, that is used for construction of this isoquant, is going through inefficient unit $C$ and spanned by axes $Oy_i$ and $Oy_j$. Unit $C$ is projected onto the weakly efficient part $AZ_k$ of the frontier. If we insert artificial unit $E$ somewhere on the ray $CD$, where point $D$ is a projection of unit $C$ on the frontier, and inside the region limited by rays $AZ_k$ and $Z_kB$, where point $B$ is a projection of point $C$ on the ray $Z_l Z_k$, then unit $C$ will be projected to the efficient point $E$ belonging to the modified frontier. Notice that components of the artificial unit $E$ coincide with corresponding components of unit $C$ except coordinates that correspond to variables $y_i$ and $y_j$. In other words, unit $E$ belongs to the section that is going through the inefficient unit $C$ and is determined by variables $y_i$ and $y_j$.

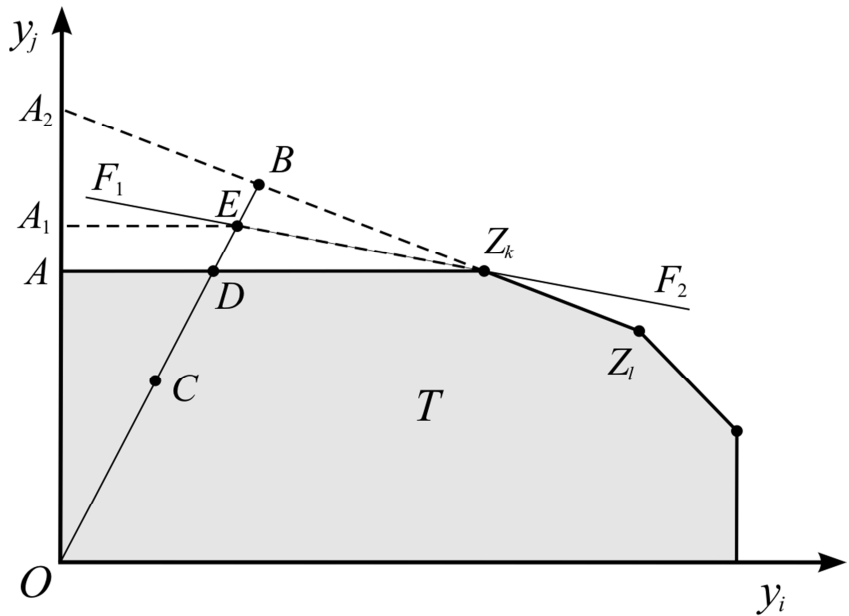

Figure 3. Removing the weakly efficient face $AZ_k$

Figure 4 shows a section of the frontier with two-dimensional plane that is going through terminal unit $Z_k$ and is spanned by axes $Ox_s$ and $Oy_p$. An artificial unit $C$ is inserted somewhere in the region limited by rays $Z_k B$ and $Z_l Z_k$. If we inserted an artificial unit $C$ somewhere in the region limited by rays $Z_l Z_k$ and $Z_k B$, then unit $Z_k$ is transformed into a non-terminal efficient unit. Again, components of artificial unit $C$ coincide with corresponding components of the terminal unit $Z_k$ except coordinates that correspond to variables $x_s$ and $y_p$. In other words, unit $C$ belongs to the section that is going through the point $Z_k$ and is determined by variables $x_s$ and $y_p$.

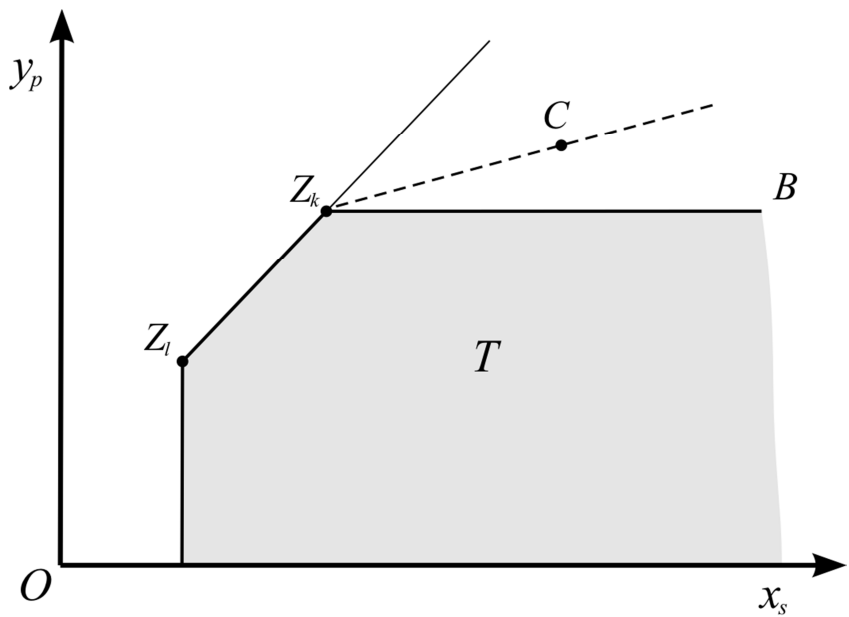

Figure 4. Section that reflects the dependence between variables $y_p$ and $x_s$

For our purpose, it is sufficient to consider only these three types of sections described above, since only these types of sections contain terminal units and terminal edges.

# 4. Algorithm for improving the frontier

Now, we describe a general scheme of the algorithm for improving the frontier in DEA models.

[INPUT:] $T_{term}$ – set of terminal units

$D_j$ – set of terminal directions for $j \in T_{term}$

Efficiency scores for all production units $j = 1, \ldots, n$

[OUTPUT:] $A$ – set of artificial units

**begin procedure** FrontierImprovement

**initialization:** $A = \varnothing$

**Part 1. Smoothing terminal units**

    **for each** $j \in T_{term}$ **do**:

        **for each** $d \in D_j$ **do**:

            **for each** two-dimensional section, that contains direction $d$ **do**:

                Insert an artificial unit on the two-dimensional section outside the PPS.

                Compute efficiency scores for all units.

                **while** the number of efficient units is less than the original number **do**:

                    Move the artificial unit closer to the frontier.

                    Compute efficiency scores for all units.

                **continue while**

            **next** two-dimensional section

            Store the new artificial unit.

        **next** $d$

    **next** $j$

Include all artificial units in the set of production units of the PPS.

**Part 2. Additional corrective steps following Part 1**

    Compute efficiency scores for all production units including artificial ones.

    Find units that were efficient and become inefficient.

    Find artificial units that caused the situations in the previous item.

    **while** there exist artificial units that have to be corrected **do**:

        Move all artificial units closer to the frontier.

        Recompute efficiency scores.

    **continue while**

    Delete inefficient artificial units.

    Compute efficiency scores for all units including also artificial production units.

**Part 3. Removing the weakly efficient faces of the frontier**

    **while** there exist units that are projected on the weakly efficient faces **do**:

        Move projection on the weakly efficient faces along the radial direction outside the PPS. Create artificial unit from such projection, and insert this artificial unit in the current iteration in the PPS.

        Compute efficiency scores.

        **while** the number of originally efficient units decreases **do**:

            Decrease the distance of the new artificial unit from the frontier.

            Recompute efficiency scores.

        **continue while**

        Store the new artificial unit.

    **continue while**

    Include all artificial units in the set of production possibility units.

**Part 4. Additional corrective steps following Part 3**

    Compute efficiency scores for all production units including artificial ones.

    Find original units that were efficient and become inefficient.

    Find artificial units that were inserted in the model for correction in the previous Part 3.

    **while** there exist artificial units that should be corrected **do**:

        Move all such artificial units simultaneously closer to the frontier along the radial direction.

        Recompute efficiency scores.

    **continue while**

    Include all efficient artificial units in the set $A$.

**end procedure**

    Observe that during the run of Part 1 and Part 3 of the algorithm some artificial units are inserted in the model. For this reason some efficient units may turn into inefficient ones since the configuration of the production possibility set (set of vertices, set of faces and their mutual disposition) may be changed. For this reason two additional stages (Part 2 and Part 4) are introduced in the algorithm in order to correct such cases. This can be accomplished by moving artificial units closer to the corresponding faces.

**Theorem.** *After the run of the algorithm the following results will be obtained:*

1) *all originally efficient units will remain efficient;*
2) *all terminal units are transformed into non-terminal efficient ones;*
3) *all inefficient units are projected onto the efficient faces of the frontier.*

Proof. Consider Part 1 of the algorithm.

In the algorithm, artificial units are generated in such a way that all efficient units (vertices) stay efficient and all terminal units turn into just efficient units. Indeed, the algorithm takes a two-dimensional section of the frontier for every terminal direction. Next, an artificial unit is inserted. Without any loss of generality, consider the two-dimensional section $S_1$ of the frontier, see Figure 2. As direction vectors of the section we took $d_1 = (e_i, 0) \in E^{m+r}$ and $d_2 = (e_j, 0) \in E^{m+r}$. Unit $Z_k$ is a vertex and ray $Z_k B$ is an edge of the polyhedral set $T$. Hence a supporting hyper-plane can be constructed in such a way that it goes through point $Z_k$ and has no other common points with set $T$. For this reason the hyper-plane cannot contain the two-dimensional section $S_1$ entirely, since in this case it would not be a supporting hyper-plane. So, this hyper-plane intersects section $S_1$ along some line $L_1 L_2$, see Figure 2. Line $L_1 L_2$ may take any position between rays $Z_k C$ and $Z_k B$. Insert an artificial unit $A$ somewhere on the ray $Z_k L_2$ between rays $Z_k B$ and $Z_k C$, no new terminal units appear except unit $A$.

According to the construction, artificial unit $A$ belongs to the two-dimensional section $S_1$ shown in the Figure 2. However some efficient units (vertices) may become inefficient. These units may be situated in the multidimensional space $E^{m+r}$. Assume that some efficient unit $F$ become inefficient after inserting artificial point $A$ in the production possibility set. Now unit $A$ is a vertex of the modified set $T$.

Let point $A_1$ be a projection of unit $A$ onto the ray $Z_k C$, see Figure 2. According to Theorem 3.15 in Preparata and Shamos (1985, p.137) unit $A$ is an extreme efficient unit and belongs to some new face $S$ of the modified set $\tilde{T}$. Unit $F$ is an interior unit of the set $\tilde{T}$. Move unit $A$ along line $AA_1$. When unit $A$ reaches point $A_1$ face $S$ turns into a face of set $T$, and unit $F$ becomes again extreme efficient. Hence there exists such point $G$ on line $AA_1$, that unit $F$ is again efficient if point $A$ belongs to the segment $GA_1$.

This implies that the algorithm can find such position on the segment $AA_1$ for a finite number of steps that all originally efficient units remain efficient and terminal units of the type $Z_k$ are transformed into non-terminal efficient units.

The algorithm processes successively all terminal units and terminal directions. Hence all terminal units are transformed into non-terminal efficient units and efficient units remain efficient.

So after the run of the first part of the algorithm all terminal units of set $T$ will turn into non-terminal efficient units (vertices) of the production possibility set, since algorithm takes all sections going through every terminal units and based on their terminal directions.

Next, consider the Part 3 of the Algorithm. Let us take section $S_2$ as an example, see Figure 3, without any loss of generality.

In this case, the directional vectors of the two-dimensional section are $d_1 = (0, e_i) \in E^{m+r}$ and $d_2 = (0, e_j) \in E^{m+r}$. Let inefficient unit $C$ be projected on the weakly efficient face, and let point $D$ be a projection of point $C$.

Take some point $E$ somewhere on the ray $CD$ between segments $Z_k A$ and $Z_k A_2$. Point $E$ is situated outside the current production possibility set $\tilde{T}$ since it lies on the ray $CD$, point $C$ is an interior point and point $D$ is a boundary point of the production possibility set (Nikaido, 1968). A supporting hyper-plane can be built in such a way that it goes through point $E$ and $Z_k$. Again, this hyper-plane intersects section $S_2$ along line $F_1 F_2$. This line may take any position between segments $Z_k A$ and $Z_k A_2$. Observe that in this case point $Z_k$ may be an intersection of several faces in the multidimensional space $E^{m+r}$ (Krivonozhko et al., 2014).

As in the previous case, there exists such position of artificial unit $E$ on the open segment $(D, B)$ that all efficient units stay efficient. This position can be found by algorithm for a finite number of steps by moving unit $E$ closer to the point $D$.

Observe that we do not need to construct some sections for this case. We use section $S_2$, Figure 3, only for explanation. In reality, artificial point $E$ can be inserted several times on segment $DB$, each time closer to point $D$, until all efficient units will remain efficient.

Such operations are repeated for every inefficient unit that is projected on the weakly efficient part of the frontier.

The case for section $S_3$ can be considered in a similar way.

Terminal unit is characterized by the fact that an infinite edge is going out from this unit. Only vectors of the following forms $d_k = (e_k, 0) \in E^{m+r}$, $k=1,\ldots,m$, $g_i = -(0, e_i) \in E^{m+r}$, $i=1,\ldots,r$, where $e_k = (0,\ldots 1,\ldots 0) \in E^m$ (the unity is in $k$th position) and $e_i = (0,\ldots,1,\ldots 0) \in E^r$ (the unity is in $i$-th position), can be the direction vectors of infinite edges of set $T$ (Krivonozhko et al. 2015, 2017). All types of sections that are used in algorithm include all such direction vectors. Hence, all types of terminal units are used in the algorithm for improving the frontier, and all types of weakly efficient faces belonging to $WEff_IT$ and/or $WEff_OT$ sets are improved by algorithm using only three types of sections. In other words, all inefficient units will be projected on the efficient parts of the frontier.

Furthermore, some efficient units may become inefficient during the solution process as a result of changing of the production possibility set configuration when some artificial units are added to the current production possibility set. For this reason Part 2 and Part 4 were included in the algorithm in order to correct such cases. New artificial units are moved closer to the frontier during the execution of Part 2 and Part 4, so all efficient units stay efficient at these stages.

This completes the proof.

We described only a general scheme of the algorithm, since it can be implemented in many different forms. For example, consider Part 1. An artificial unit $A$ is inserted somewhere between rays $Z_k B$ and $Z_k C$ outside the current set $T$, see Figure 2. So, we can implement in the algorithm the following different cases.

1. To construct ray $Z_k B$ explicitly as a continuation of ray $Z_l Z_k$, see Figure 2.
2. To place unit $A$ somewhere between ray $Z_k C$ and axis $Ox_i$.

The first case requires additional computations in order to find segment $Z_l Z_k$. In the second case, we can place unit $A$ by chance outside the region limited by rays $Z_k B$ and $Z_k C$. In this case some initially efficient units may become inefficient ones, then distance $AA_1$ is decreased and computations are repeated.

In our implementation of the algorithm we use the second case. Computational results are given in the next section.

Moreover, experts can be used in order to form artificial units. The version of the algorithm, which was exposed in the paper, is deliberately described in the form that works automated. These were specifically done, because participation of decision makers or experts may increase the computational time and insert some slips in frontier improvement. However the algorithm enables one to involve experts or decision making persons in the process of construction of the new frontier. In fact, the algorithm generates only main "frame" for frontier improvement. This frame may be used to generate artificial units with necessary operating mixes and to place them in a specific location. Consider Figure 2, after the work of the algorithm experts can insert instead of unit $A$ an artificial unit in the cone determined by rays $Z_k A$ and $Z_k A_1$, where point $A$ is the last position of this point after the work of the algorithm. The similar operations can be done for output isoquants, see Figure 3, and for $x_s$-$y_p$ sections, see Figure 4. The frame of the algorithm guarantees that the efficient frontier will not be destroyed.

# 5. Computational experiments

In our computational experiments we used the software FrontierVision, a specifically developed program by our team for the DEA models (Krivonozhko et al. 2014). This program allows us to visualize the multidimensional production possibility set by means of constructing two- and three-dimensional sections of the frontier.

At first, we took data from 174 Russia bank's financial accounts for January 2009. We used the following variables as inputs: working assets, time liabilities, and demand liabilities. As output variables we took: equity capital, liquid assets, fixed assets. Max, min and mean statistics are presented in Table 1.

Figure 5 represents two input isoquants for unit 149 that are intersections of the six-dimensional production possibility set with two-dimensional planes for unit 149. This unit in the figure is shown by white circle. Other small circles represent orthogonal projections of actual and artificial units onto the section. The red color means that the corresponding unit is efficient. The green and yellow colors denote units with low and intermediate values of efficiency score. The curve 1 shows input isoquant for the original set $T$. The curve 2 is built for the transformed set $T$. Directions of the two-dimensional plane are determined by the following inputs: demand liabilities and working assets.

Table 1. Data for Russian banks, January 2009

| Variables | Mean | St. deviation | Min | Max |
|---|---|---|---|---|
| *Outputs* | | | | |
| $Y_1$ – Liquid assets, bln roubles | 6.62 | 9.36 | 0.18 | 63.9 |
| $Y_2$ – Equity capital, bln roubles | 3.52 | 4.82 | 0.52 | 26.9 |
| $Y_3$ – Fixed assets, bln roubles | 1.02 | 1.32 | 0.01 | 6.9 |
| *Inputs* | | | | |
| $X_1$ – Demand liabilities, bln roubles | 15.16 | 18.16 | 0.55 | 105.0 |
| $X_2$ – Time liabilities, bln roubles | 27.06 | 37.86 | 0.34 | 191.6 |
| $X_3$ – Working assets, bln roubles | 36.49 | 47.34 | 1.83 | 249.2 |

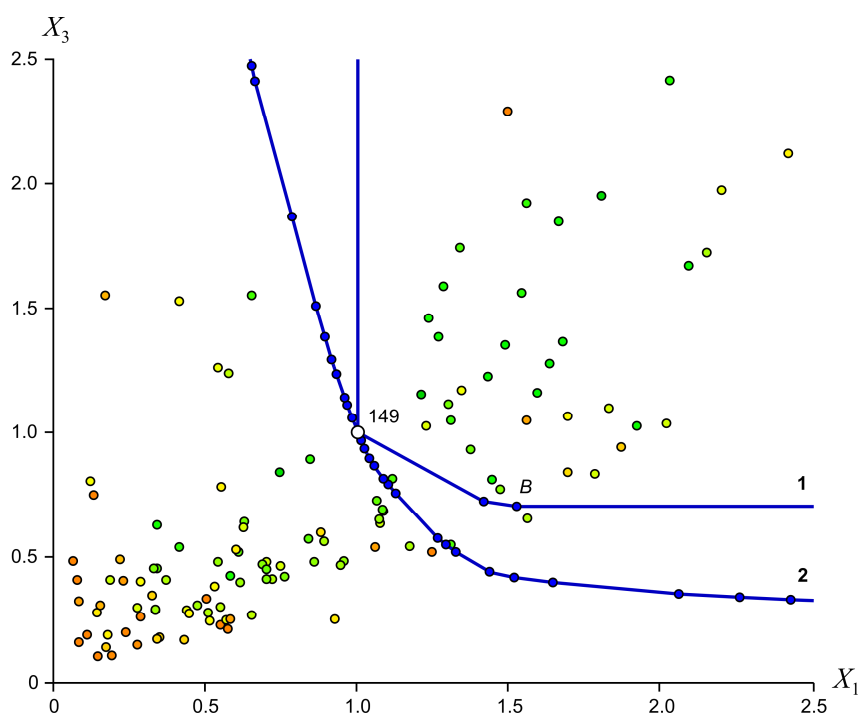

Figure 5. Input isoquant for unit 149

Input isoquants are constructed for efficient unit 149. However, they have only one common point, unit 149, this means that curve 1 consists mainly of weakly efficient points of the frontier except unit 149. It seems to us that some units in the figure are situated outside the production possibility set. Certainly this is not so, because the figure shows orthogonal projections of other production units being in the six-dimensional space onto the two-dimensional plane of section. This feature of our software allows experts to estimate dispositions of units in the multidimensional space of inputs and outputs. This feature can be used by experts optionally.

Figure 6 depicts two output isoquants for unit 53, these curves are intersections of the six-dimensional production possibility set with two-dimensional plane for unit 53. Directions of the plane are taken as follows: fixed assets and liquid assets. The curve 1 shows input isoquant for set $T$ and curve 2 is input isoquant for the transformed set $T$. It is worthy of note that all points of the curve 1 are weakly efficient points, this means that all efficiency scores that are measured relative to these points would be given incorrect values.

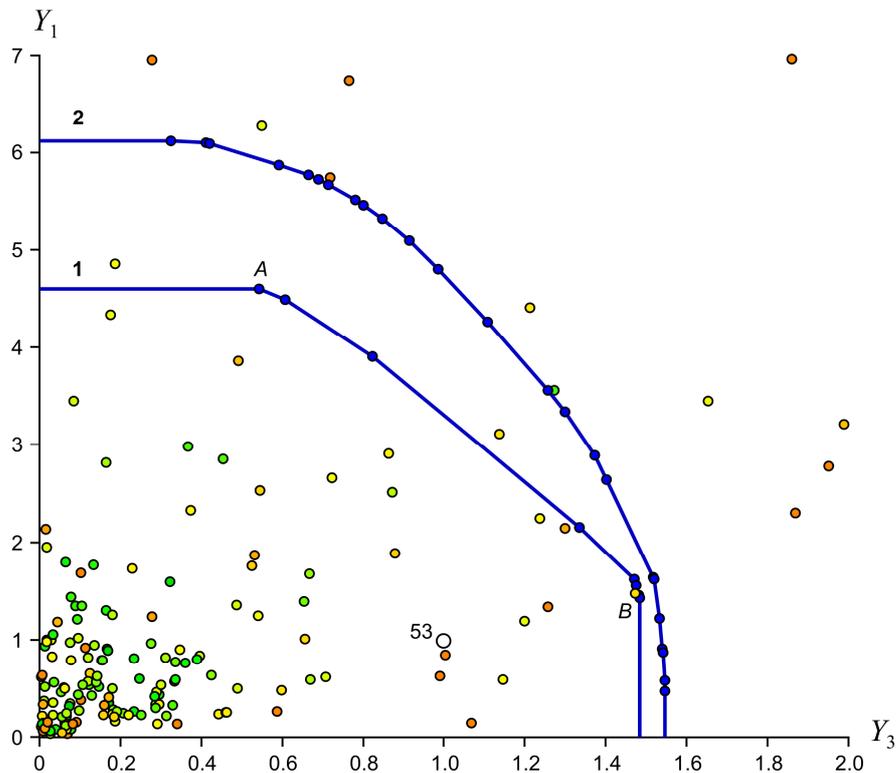

Figure 6. Output isoquant for unit 53

The number of originally efficient units among banks is 26. The number of inefficient units is 148. Almost all inefficient units are projected onto the weakly efficient parts of the frontier, set $WEffT$ and/or $WEff_OT$, this number is equal to 146.

The algorithm inserted 412 artificial units in the original set of units in order to improve the frontier. The algorithm takes approximately 10 minutes to improve the frontier on the personal computer with Intel Core i3 CPU 3.33 GHz.

Output isoquants are built for inefficient unit 53. Curves 1 and 2 have no common points. This implies that curve 1 originates from weakly efficient points of the original frontier only. So, curve 2 is completely outside curve 1.

After running the algorithm, all inefficient units are projected onto efficient parts of the frontier. Bougnol and Dulá (2009) observed that almost all extreme efficient units are anchor units. So, it would be interesting to check: how many inefficient units are projected onto the weakly efficient faces of the frontier in other real-life models.

To achieve this purpose, we expanded our computational experiments and included in it two additional datasets. We took the data for electricity utilities on Sweden 1987, see (Førsund et al., 2007). The number of production units in this model is 163, among these units there are 110 inefficient units. And 109 inefficient units are projected on the weakly efficient parts of the frontier.

As the second additional dataset, we took data of the nursing and home care sector of Norwegian municipalities. There are three inputs and ten outputs in this model, see details in Krivonozhko et al. (2014). This model has 469 original units. Computations on the BCC model show that there are 129 efficient units among them and 340 inefficient ones. All these inefficient units are projected onto the weakly efficient faces of the frontier. Hence the frontier in the DEA models is really in need of improving.

## 6. Conclusions

In this paper, we proposed an algorithm for improving the frontier in the DEA models. Strictly speaking, only a general scheme of the algorithm was described in detail, since it can be implemented in many different forms that depend on original program modules used for its construction.

The computational results after improving the frontier are more accurate and reliable, because efficiency scores of inefficient units are determined now relative to efficient points on the frontier. Moreover, some marginal rates at weakly efficient points on the frontier may take infinite value or zero value before improving the frontier, what is unnatural for economic processes. Figures 5 and 6 in the paper show such situations. In Figure 5, the left marginal rate of substitution at point 149 belonging to curve 1 is equal to infinity, the right marginal rate of substitution at point $B$ is equal to zero. We can see the similar situation for marginal rate of transformation. The left marginal rate of transformation at point $A$ on the curve 1, see Figure 6, is equal to zero. The right marginal rate of transformation at point $B$ on the curve 1 is equal to infinity, see Figure 6. After improving the frontier marginal rates take finite values at points on the frontier, which envelop the observed production units. Observe that it is impossible to remove all terminal units in DEA models. The algorithm inserts new terminal units farther from the origin of coordinates, for this reason all inefficient observed production units are projected on the efficient parts of the frontier, as it supposed in the DEA approach. Figures 5 and 6 illustrates this fact clearly.

The point is that the number of observed units is always finite, and this number is typically not sufficient for construction of the complete unknown production possibility set. Artificial units are introduced in order to transform weakly efficient faces into efficient ones. However, all originally efficient faces will remain efficient. This is the essence of the proposed algorithm.

Economic functions are nice instruments for decision making in economics and business. However, they are seldom used based on DEA models. We developed the software FrontierVision that enables us to construct various economic functions in DEA models as sections of the frontier with two-dimensional subspaces (Krivonozhko et al. 2004). It turned out that a large share of such functions look angular and awkward, see for example Figure 5 and Figure 6 in the paper. This is because the number of efficient observed units is not sufficient for construction of the production possibility set, as was explained above. So, another benefit of improving the frontier is that it is more suitable now as instruments for decision making.

Computational experiments using real-life datasets confirmed that the algorithm works reliably and improves the frontier significantly.

## Acknowledgement

This work was supported by the Russian Science Foundation (project No. 17-11-01353).